\documentclass[english,draft]{article}
\usepackage{amssymb,amscd,amsthm,amsmath,babel,amsfonts,makeidx,xypic}

\theoremstyle{plain}
\newtheorem{theorem}{Theorem}[section]
\newtheorem{lemma}[theorem]{Lemma}
\newtheorem{corollary}[theorem]{Corollary}
\newtheorem{proposition}[theorem]{Proposition}

\theoremstyle{definition}

\newtheorem{definition}[theorem]{Definition}
\newtheorem{remark}[theorem]{Remark}

\newtheorem{example}[theorem]{Example}

\numberwithin{equation}{section}

\newcommand\opn[2]{%
    \newcommand{#1}{\operatorname{#2}}}

\newcommand\NN{{\mathbb N}}

\newcommand\PP{{\mathbb P}}

\newcommand\frk{\mathfrak}

\newcommand\mm{{\frk m}}

\newcommand\qq{{\frk q}}

\opn\charac{char}
\opn\projdim{proj\,dim}
\opn\depth{depth}
\opn\rank{rank}
\opn\rlex{rlex}
\opn\fine{end}
\opn\supp{supp}
\opn\Ass{Ass}
\opn\Proj{Proj}
\opn\Spec{Spec}
\opn\Soc{Soc}
\opn\reg{reg}
\opn\Md{Md}
\opn\dBorel{dBorel}
\opn\Borel{Borel}
\opn\Shad{Shad}
\opn\Span{span}
\opn\Hs{Hilb}
\opn\ini{in}
\opn\Gin{Gin}
\opn\grade{grade}
\opn\greg{greg}
\opn\Ht{ht}

\opn\Ker{Ker}
\opn\coker{coker}

\opn\im{Im}
\opn\Hom{Hom}
\opn\Tor{Tor}
\opn\Ext{Ext}
\opn\id{id}
\opn\Id{Id}
\opn\Homst{^*Hom}
\opn\Extst{^*Ext}
\opn\Gamst{^*\Gamma}
\opn\Hst{^*H}

\opn\lk{lk}
\opn\st{st}

\let\:=\colon
\let\phi=\varphi

\newcommand\pnt{{\raise0.5mm\hbox{\large\bf.}}}

\newcommand\ra{\rightarrow}

\newcommand\sat{^{sat}}
\opn\po{pol}
\opn\pol{^{\bf p}} 

\begin{document}
\noindent  
{\Large Characteristic-free bounds for the Castelnuovo-Mumford regularity}\\

\vspace{.3cm}
\noindent
{\bf Giulio Caviglia}\footnote{Partially supported by I.N.d.A.M., Rome}. (\texttt{caviglia@math.ukans.edu})\\
{\footnotesize Department of Mathematics, University of Kansas, Lawrence, KS 66045, USA}\\
\noindent
{\bf Enrico Sbarra}. (\texttt{sbarra@dsm.univ.trieste.it})\\
{\footnotesize DSM, Universit\`a di Trieste, Via A. Valerio 12/1, I - 34127,  Trieste}\\

\noindent
{\textsl \today}\\\\

\noindent
{\small {\bf Abstract}\\
We study bounds for the Castelnuovo-Mumford regularity of homogeneous ideals in a polynomial ring in terms of the number of variables and the degree of
the generators. In particular our aim  is to give a positive answer to a question posed by  Bayer and Mumford in \cite{BaM}, by 
showing that the known upper bound in characteristic zero holds true also in positive characteristic. We first analyze Giusti's proof, 
which provides the result in characteristic $0$, giving some insight on the 
combinatorial properties needed in that context. 
For the general case  we provide a new argument which employs Bayer and Stillman criterion for detecting regularity.}\\\\

\section*{Introduction}
The Castelnuovo-Mumford regularity is an important invariant in commutative algebra and algebraic geometry, which gives an estimate of 
the complexity of 
computing a minimal free resolution. It is common in the literature to attempt to find bounds for this invariant and, in general, the expected results range 
quite widely, from the well-behaved examples coming from the algebraic geometry, as suggested by the Eisenbud-Goto Conjecture \cite{EiGo}, to the worst case 
provided by the  example of Mayr and Meyer \cite{MaMe}. Clearly, when the assumptions are quite unrestrictive, the regularity can be very large. If one works
with an homogeneous ideal $I$ in a polynomial ring $R=K[X_1,\dots,X_n]$ over a field $K$, a very natural question is to ask whether the regularity 
can be limited just by knowing that the ideal is generated in degree less than or equal to some positive integer $d$. What was known to this point were bounds 
depending on the characteristic of the base field $K$. If $\charac K=0$, as observed in \cite{BaM}, Proposition 3.8, from the work of Giusti \cite{Gi} and
and Galligo \cite{G}, \cite{G1} one can derive
\begin{equation}\tag{A}\label{char0}
\reg(I)\leq (2d)^{2^{n-2}},
\end{equation}
which seems to be sharp (see again \cite{MaMe}).

On the other hand, in any characteristic, it has been proven in \cite{BaM}  using ``straightforward cohomological methods'' that 

\begin{equation}\tag{B}\label{charp}
\reg(I)\leq (2d)^{(n-1)!},
\end{equation}
but in the same paper it is asked whether \eqref{char0} holds in general independently of the characteristic. The main purpose of this article is to give a 
positive answer to this question.
The main effort in extending the result to positive characteristic is that this proof utilizes the combinatorial structure of the 
{\em generic initial ideal} in characteristic zero.

The generic initial ideal has been introduced in \cite{G1}, where it was defined  with the assumption that the base field has 
characteristic zero. The definition has been generalized a few years later for base fields of any characteristic in \cite{BaSt}, and grew in
importance, as many recent results demonstrate. 
One of the points of major interest in considering the generic initial ideal $\Gin_{\rlex}(I)$ with respect to the (degree) reverse lexicographic order 
of an homogeneous ideal $I$ is that this is a monomial ideal with the same Hilbert function, projective dimension  and Castelnuovo-Mumford 
regularity as $I$. Furthermore a generic initial ideal is Borel-fixed, i.e. it is invariant under the action of the Borel group, which is the subgroup of 
$GL_n(K)$ consisting of all non-singular upper-triangular $n\times n$  matrices with coefficients in $K$. 
According to the characteristic of the underlying field, whether it is $0$ or positive, Borel-fixed ideals have a more or less manageable combinatorial structure.
It may be convenient to recall some of the most interesting notions used in this context 
and  to fix some terminology since this is not unique in the literature.\\
 We refer the interested reader to the detailed treatise in \cite{P} and \cite{P1} 
for further information.\\
Let $K$ be an infinite field (which is not a restrictive hypothesis for our purposes). 
Given a monomial $u$ we denote $\max\{i\:X_i\mid u\}$ by $m(u)$.
Let now $p$ be a prime number and $k$ an integer. The {\it $p$-adic expansion} of $k$ is the expression of $k$ as $\sum_ik_ip^i$, with $0\leq k_i \leq p-1$. If $k=\sum_ik_ip^i$ and
$l=\sum_il_ip^i$ are the $p$-adic expansions of the two integers $k$ and $l$ respectively, one sets $k\leq_p l$ iff $k_i\leq l_i$ for all $i$.\\
First of all notice that an ideal $I$ which is fixed under the action of the Borel group (i.e. a {\it Borel}({\it-fixed}) ideal) is monomial.\\
A {\it standard Borel}({\it-fixed}) 
(or {\it strongly stable}) ideal $I$ is an ideal endowed with the following property: for every $u\in I$, if $X_i\mid u$ then  $\frac{X_ju}{X_i}\in I$, 
for every  $j<i$.\\
It is wider the class of ${\it stable}$ ideals defined by a weaker exchange condition on the variables of the monomials: An ideal $I$ is stable 
iff for every $u\in I$, $\frac{X_ju}{X_{m(u)}}\in I$, for every  $j<m(u)$.\\
Finally, an ideal $I$ is said to be {\it p-Borel} iff
for every monomial $u\in I$, if $l$ is the maximum integer such that $X_i^l|u$, then $\frac{X_j^ku}{X_i^k}\in I$, for every $j<i$ and $k\leq_p l$.\\
Standard Borel ideals are Borel. If $\charac K=0$ every Borel ideal is standard. If $\charac K=p$ a monomial ideal is Borel iff it is p-Borel.\\
\noindent
The crucial difference between characteristic $0$ and positive characteristic can be probably noticed at first glance and consists of the fact 
that the combinatorial structure of standard Borel ideals is easier than that of the non-standard ones, which relies on the $p$-adic expansion of non-negative 
integers. This difference in behaviour also results in the fact that
there is an complete description of the minimal free graded resolution of a standard Borel ideal $I$ in terms of the monomials of its minimal system of 
generators $G(I)$ (cf. \cite{EK}, \cite{AH}), while the task of finding an analogous for non-standard Borel ideals still seems to be too difficult.
In particular, the graded Betti numbers of a standard Borel ideal can be computed explicitly in terms of $G(I)$ and it is  easily deduced that its
Castelnuovo-Mumford regularity equals the highest degree of an element in $G(I)$, i.e. the so-called {\em generating degree} of $I$.\\
Thus, it is quite clear that the assumption of characteristic $0$ has made the task of investigating the generic initial ideal easier.\\

\noindent
This paper is organized in two sections. The first section is dedicated to a better understanding of the arguments that lead to the inequality expressed in
\eqref{char0}, which are due to Galligo \cite{G} and Giusti \cite{Gi}. This resulted in the introduction and enquiry of certain ideals, which we call 
{\em weakly stable}. As a result we obtain a bound for the regularity which improves \eqref{char0}. 
In the second section we capitalize on these ideas and we prove 
a formula that relates the regularity of the beginning  ideal with its
generating degree and the regularity of its sections by an  
almost-regular sequence of linear forms (Theorem \ref{main}).
As a consequence, we obtain the desired result (Corollary \ref{main2}).

\section{Weakly stable ideals}
The purpose of this section is to study a certain class of monomial ideals which we call {\it weakly stable}.\\ 
It may be appropriate before proceeding to recall the definition of the Castelnuo\-vo-Mumford regularity, 
whereas we refer the reader to \cite{EiGo}, \cite{Ei} and \cite{BS} for further details on the subject.

\begin{definition}
Let $M$ be a finitely generated graded $R$-module and let $\beta_{ij}(M)$  denote the graded Betti numbers of $M$ (i.e. the numbers
$\dim_K\Tor_i(M,K)_j$). The  {\it Castelnuovo-Mumford regularity} $\reg(M)$ of $M$ is $$\max_{i,j}\{j-i\:\beta_{ij}(M)\neq 0\}.$$
\end{definition}
\noindent
Remember also this useful characterization of regularity in terms of the local cohomology modules of $M$, which we shall  use in the following. 
Since the graded local cohomology modules $H^i_\mm(M)$ with support in the maximal graded ideal $\mm$ of $R$ are Artinian,
one defines $\fine(H^i_\mm(M))$  to be the
maximum integer $k$ such that $H_\mm^i(M)_k\neq 0$. Thus
$$\reg(M)=\max_i\{\fine(H^i_\mm(M))+i\}.$$
Finally, a finitely generated $R$-module $M$ is said to be  {\em $m$-regular} for some integer $m$ iff $\reg(M)\leq m$.\\
Let us now underline which are the main steps necessary for the proof of \eqref{char0}. 
For any field $K$ and any homogeneous ideal $I$ it is well-known that $\reg(I)$ equals $\reg(\Gin_{\rlex}(I))$, and if $\charac K=0$ then $\Gin_{\prec}(I)$ is a 
standard Borel ideal for any term order $\prec$ for which $X_1\succ X_2\succ\ldots\succ X_n$, so that its regularity equals its generating degree $D$. 
Furthermore, if $\charac K=0$, the so-called Crystallisation Principle (which from now on we shall abbreviate with {\bf CP}) holds:
\begin{itemize}\item[{\bf CP}:]
Let $I$ be an homogeneous ideal generated in degrees $\leq d$. Assume that $\Gin_{\rlex}(I)$ has no generator in degree $d+1$.  
Then $D\leq d$ (cf. \cite{Gr}, Proposition 2.28).
\end{itemize}
Hence, thanks to the good properties of $\Gin_{\rlex}$ and an induction argument on the numbers of variables,
one obtains bounds for $D$ in terms of the generating degree $d$ of $I$ (cf. \cite{Gi}, in particular 
the ``Proof of Theorem B''), and this completes the argument.\\
One notices that in the proof the hypothesis $\charac K=0$ is used solely for the purpose of  exploiting the combinatorial structure of $\Gin_{\rlex}$. \\
Furthermore {\bf CP} only holds true in characteristic $0$. Consider for instance the ideal $(X^{2p},Y^{2p})$ in $K[X,Y]$ with $\charac K=p\neq 2$. 
Then $\Gin_{\rlex}(I)=(X^{2p},X^pY^p,Y^{3p})$. Here it is sufficient  to observe that the ideal $(X^{2p},Y^{2p})$ is the ideal 
generated by the images of $X^2$ and $Y^2$ under the Frobenius map $R \ra R$, $X\ra X^p$. In fact the following more general result is well known.

\begin{proposition}
Let $I$ be an homogeneous ideal of $R=K[X_1,\ldots,X_n]$ with $\charac K=p$ and let $F$ be the Frobenius map. Then for any term order $\tau$ one has
$$\Gin_{\tau}(F(I))=F(\Gin_{\tau}(I)).$$
\end{proposition}
\begin{proof}
Note that the computation of the initial ideal of $F(I)$ can be performed in $K[X_1^p,\ldots,X_n^p]$, i.e. the S-pairs of $F(I)$ are just the $p$-th
power of the S-pairs of $I$, so that $F(\ini_{\tau}(I))=\ini_{\tau}(F(I))$. This suffices, since by definition $\Gin_{\tau}(F(I))=\ini_{\tau}(g(F(I)))
=\ini_{\tau}(F(g(I)))$, where $g$ is a generic change of coordinates.
\end{proof}

Since {\bf CP} fails in positive characteristic, one might wonder if the proof of \eqref{char0} could be 
performed by making use  of {\it lexicographic} (also called {\it lex-segment}) ideals instead of generic 
initial ones. In fact 
there is a natural counterpart of {\bf CP} that is expressed by the so-called Gotzmann's Persistence Theorem (see \cite{Gz} or \cite{Gr}, Theorem 3.8): 
Given an ideal $I$ with generating degree $D$ the lexicographic ideal $L$ associated with $I$ cannot have generators
in degree $k>h$ for any $k$ if it has none in degree $h>D$.\\
 This property though is not strong enough to conclude the proof of the bounds since
there is something missing: Modding out the last variable the resulting lexicographic ideal in the smaller polynomial ring does not 
fulfill the 
Gotzmann's Persistence Theorem for the same $D$, so that if one would
like to proceed recursively, one would obtain much higher bounds.
This problem led to the understanding of the property  one needs to recover Giusti's argument, and motivates the following definition.\\
As before we let $m(u)$ denote  $\max\{i\:X_i\mid u\}$ for any monomial $u$.

\begin{definition}
A monomial ideal $I$ is called {\em weakly stable} iff the following property holds.
For all $u\in I$ and for all $j< m(u)$ there exists a positive integer $k$ such that 
$\frac{X_j^ku}{X_{m(u)}^l}\in I$, where $l$ is the maximum integer such that $X_{m(u)}^l\mid u$. 
\end{definition}

\noindent
It is straightforward from the definition that if $I$ is weakly stable so is $\bar{I}$, the quotient ideal 
obtained  by modding out the last variable (any other variable would do after re-labelling). One 
easily verifies that finite intersections, sums and products of weakly stable ideals are weakly stable.
It is worth to point out that a monomial ideal is weakly stable iff its associated prime ideals are 
lexicographic, i.e. of the form $(X_1,\ldots,X_i)$ for some $i$. This and other combinatorial properties 
of weakly stable ideals have been proven in \cite{C}. We also recall that monomial ideals which 
have lexicographic associated prime ideals have been used in \cite{BeGm} for algorithmic computations of the
Castelnuovo-Mumford regularity.
 
\begin{remark}\label{evvabe}
Strongly stable, stable and $p$-Borel ideals are weakly stable.
\end{remark}

Henceforth we let $D(I)$ denote the generating degree of the ideal $I$, i.e. the maximum of the degrees 
of a minimal set of generators of $I$.
We next prove a bound on the cardinality of the minimal system of generators of a weakly stable ideal in terms of the generating degrees 
of its sections.
In the following we shall denote by $I_{[i]}$ the image of an  ideal $I$ in  $R/(X_{i+1},\ldots,X_n)$, for $i=1,\ldots,n-1$ and we let
 $I_{[n]}=I$.\\
Given a monomial $u$ in $K[X_1,\ldots,X_n]$ we denote the maximal degree of the variable $X_i$ appearing in $u$ by 
$\Md_i(u)\doteq\max\{j\:X_i^j\mid u\}$,
for $i=1,\ldots,n$. Accordingly, if $J$ is a monomial ideal  we define $\Md_i(J)\doteq\max_{u\in G(J)}\{j\:X_i^j\mid u\}$ which
is $\max_{u\in G(J)}\{\Md_i(u)\}$, for $i=1,\ldots,n-1$.
For the next proposition it is crucial what  we prove in the next lemma.

\begin{lemma}
Let $I$ is a weakly stable ideal. Then  $\Md_i(I_{[i]})=\Md_i(I)$ for all $i=1,\ldots,n-1$.
\end{lemma}
\begin{proof}
It is immediate that $\Md_i(I_{[i]})\leq\Md_i(I)$. Suppose now that 
$\Md_i(I_{[i]})<\Md_i(I)$ and let us find a contradiction. For the
sake of simplicity let $s\doteq\Md_i(I_{[i]})$. Then there exists 
$u\in G(I)$ such that 
$X_i^{s+1}\mid u$ and $m(u)>i$. 
Let us choose such a counterexample such that $m(u)$ is the smallest possible. It thus follows that there exists an integer $k$ such 
that $\frac{uX_i^k}{X_{m(u)}^k}$ is an element of $I$, because $I$ is weakly stable.
Hence there exists $n\in G(I)$ such that $n\mid \frac{uX_i^k}{X_{m(u)}^k}$ and $m(n)<m(u)$. Therefore $\Md_i(n)\leq s$, 
so that $n\mid uX_i^k$ and $n\nmid u$. 
But this implies that $s\geq \Md_i(n)\geq\Md_i(u)+1$ which is $\geq
s+2$ and we are done.
\end{proof}

\begin{proposition}\label{weakly1}
Let $I\subset R=K[X_1,\ldots,X_n]$, with $n\geq 2$, be a weakly stable ideal. Then $$|G(I)|\leq \prod_{i=1}^{n-1}(D(I_{[i]})+1).$$
\end{proposition}
\begin{proof}
Keeping in mind that we can establish whether two monomials of $G(I)$ are distinct just by looking at their first $n-1$ variables,
it is clear that $|G(I)|\leq \prod_{i=1}^{n-1}(\Md_i(I)+1)$. We have already verified that, if $I$ is weakly stable, then
 $\Md_i(I)=\Md_i(I_{[i]})$ for all $i=1,\ldots,n-1$. This is sufficient because $\Md_i(I_{[i]})$ is obviously $\leq D(I_{[i]})$
for all $i=1,\ldots,n-1$.
\end{proof}

\noindent
In the following we consider weakly stable ideals which fulfill the following condition {\em with respect to $d$}:
\begin{itemize}\item[(*)]
There exists an integer $d$ such that, for all $i \geq d$, 
if $I$ has no minimal generator of degree $i$ then it also has none of 
degree $i+1$.
\end{itemize}

\begin{proposition}\label{troubled}
Let $I\subset R=K[X_1,\ldots,X_n]$, with $n \geq 2$,  be a weakly stable ideal for which Condition (*) holds w.r.t. $d$. Then 
$$D(I)\leq d-1+\prod_{i=1}^{n-1}(D(I_{[i]})+1).$$
\end{proposition}
\begin{proof}
The hypothesis on the generators implies that $D(I)-d+1$ is less than or equal to the cardinality of $G(I)$ and since $I$ is weakly 
stable the assertion follows directly from Proposition \ref{weakly1}.
\end{proof}
Let us now look at a smaller class of weakly stable ideals, for which Condition (*) holds w.r.t. $d$ and such that the quotient ideals 
verify Condition (*)
w.r.t. the same $d$. In other words we consider weakly stable ideals which verify the following condition:
\begin{itemize}\item[(**)]
$I_{[i]}$ verifies Condition (*) w.r.t. $d$ for all $1\leq i\leq n$. 
\end{itemize}

\noindent
Such ideals exist and as main example we consider the generic initial ideal $\Gin_{\rlex}(I)$ of an homogeneous 
ideal $I$ in $K[X_1,\ldots,X_n]$, where the characteristic of $K$ is $0$
and the generating degree of $I$ is less than or equal to $d$. 
As noticed already many a time, $\Gin_{\rlex}(I)$ is strongly stable and {\it a fortiori} weakly stable as observed in Remark 
\ref{evvabe}. Condition
(*) is verified for such an ideal by virtue of {\bf CP}, whereas Condition (**) holds because, since the chosen term order is the 
reverse lexicographic order, 
one has
$\Gin_{\rlex}(I)_{[i]}=\Gin_{\rlex}(I)+(X_{i+1},\ldots,X_n)$ which is the same as 
$\Gin_{\rlex}(I+(X_{i+1},\ldots,X_n))=\Gin_{\rlex}(I_{[i]})$, and {\bf CP} applies since the generating degree of $I_{[i]}$ is obviously 
$\leq d$.

\begin{corollary}\label{CC}
Let $I\subseteq K[X_1,\ldots,X_n]$, with $n\geq 2$,  be a weakly stable ideal, which fulfills Condition (**) w.r.t. $d$. Then
$$D(I)\leq (2d)^{2^{n-2}}.$$
\end{corollary}
\begin{proof}
According to Proposition \ref{troubled}, one has  that $D(I_{[i]})\leq B_i$, where we let $B_1\doteq d$ and $B_i\doteq d-1+\prod_{j=1}^{i-1}(B_j+1)$,
for all $i>1$. One can easily verify that $B_2=2d$ and 
also determine any element of the sequence by its previous one by the following easy computation
$B_i=(B_{i-1}-(d-1))(B_{i-1}+1)+d-1=B_{i-1}^2-(d-2)B_{i-1}$. Therefore, since we may assume without any loss of generality that $d\geq 2$, we have
$B_i\leq B_{i-1}^2$. Thus, for all $i\geq 2$, we have that  $B_i\leq (2d)^{2^{i-2}}$. In particular $D(I)\leq B_n\leq (2d)^{2^{n-2}}$ and this proves the
statement of the corollary.
\end{proof}
\noindent
The bound for the regularity expressed in \eqref{char0} follows now easily under the assumption $\charac K=0$. 

\begin{corollary}\label{main3}
Let $I$ be an ideal of $R=K[X_1,\ldots,X_n]$ with $n\geq 2$ and $\charac K=0$. Let $I$ be generated in degree $\leq d$. Then 
$$\reg(I)\leq (2d)^{2^{n-2}}.$$
\end{corollary}
\begin{proof}
Recall that $I$ and $\Gin_{\rlex}(I)$ have the same regularity and that the latter is a stable ideal, therefore its regularity equals its generating degree. 
By the observations preceding Corollary \ref{CC}, this ideal fulfills the hypotheses of the previous corollary  
and the conclusion results from its straightforward application.
\end{proof}

\section{Bounds for the regularity}
In this section we show that the well-known doubly exponential bounds hold independently of the characteristic. This improves \cite{BaM}, Theorem 3.7 and 
Proposition 3.8.
The techniques used here are based upon general properties of local cohomology and almost-regular sequences of linear forms.\\
Henceforth, by Flat Extension, we may assume without loss of generality that $|K|=\infty$.\\
We notice here a few easy facts which are used in the rest of the section. Given an arbitrary ideal $I$ we let $I\sat$ denote 
the {\em saturation}
$I:\mm^\infty=\cup_{k\geq 0}I:\mm^k$ of $I$ with respect to $\mm$. Let $I,J$ be two arbitrary ideals. If $I\subseteq J\subseteq I\sat$ 
then $J\sat=I\sat$. Recall that, given a finitely generated $R$-module
$M$ an homogeneous element $l\in R_d$ is said to be  
{\em almost-regular for $M$} iff the multiplication map $\diagram
M_k\rto^{\cdot l} & M_{k+d} \enddiagram$ is injective for all $k\gg 0$. 
We say that $l_1,\dots,l_r$ form  an {\em almost-regular sequence for
$M$} if $l_1$ is almost-regular for $M$ and $l_{i+1}$ is almost-regular for $M/(l_1,\ldots,l_i)M$ for all $i=1,\ldots,r-1$.
One can show that a homogeneous form is almost regular for an $R$-module $M$ iff it is not
contained in any associated prime ideal of $M$ other than the homogeneous maximal ideal.
Recall also that, since we are assuming that the cardinality of $K$ is $\infty$, if $M$ has positive
dimension then any generic form is almost-regular.   

\begin{remark}\label{somehelp}
Suppose that $l$ is an almost-regular form for $R/I$. Then $I:l^\infty=I\sat$.
In fact, let us consider $I=\qq_1\cap\qq_2\cap\ldots\cap\qq_{r-1}\cap\qq_r$ a minimal primary decomposition of $I$, 
where $\qq_r$ is the component associated with the maximal
ideal. Since $l$ is almost-regular one has that $l\not\in\cup\sqrt{\qq_i}$, 
where the union is taken over all indices $i<r$. Thus 
$I\sat=\cap_{i=1}^r\qq_i:\mm^\infty=\cap_{i=1}^{r-1}(\qq_i:\mm^\infty)$. On the other hand, it is enough to show that 
$I:l^\infty\subseteq I\sat$. But $I:l^\infty=\cap_{i=1}^r\qq_i:l^\infty=\cap_{i=1}^{r-1}(\qq_i:l^\infty)=\cap_{i=1}^{r-1}\qq_i$.\\
{\em A fortiori} we also have $I:l^k\subseteq I\sat$ for all $k\in\NN$.
\end{remark}

Before introducing the main result of this section we prove two useful lemmata.
In the following we let $\lambda(\cdot)$ denote the function length.

\begin{lemma}\label{snakesnake}
Let $l$ be an element of $R$ and  $I$ an ideal of $R$. For any integer $a\geq 0$ one has
$$\lambda\left(\frac{I:l^a}{I:l^{a-1}}\right)=\lambda\left(\frac{I:l^a+(l)}{I:l^{a-1}+(l)}\right)
+\lambda\left(\frac{I:l^{a+1}}{I:l^a}\right),$$
whenever all of the above lengths are finite.
\end{lemma}
\begin{proof}
Consider the following exact sequence
$$\diagram 0 \rto & \dfrac{I:l^{a+1}}{I:l^a}  \rto^{\cdot l} &
\dfrac{I:l^a}{I:l^{a-1}}   \rto & \dfrac{I:l^{a}}{(I:l^{a-1})+l(I:l^{a+1})}
\rto & 0\enddiagram,
$$ where the third term is 
$$ \frac{I:l^{a}}{(I:l^{a-1})+(l)\cap (I:l^{a})}\simeq \frac{I:l^a+(l)}{I:l^{a-1}+(l)}.$$
The conclusion follows immediately from the additivity of the function length.
\end{proof}

\begin{lemma}\label{scianti}
Let $I$ be an homogeneous ideal and $l$ be an almost-regular linear form. Let $k$ be the smallest integer such that $I:l^\infty=I:l^k$. Then
$k\leq\reg(I)$.
\end{lemma}
\begin{proof}
Clearly it is enough to prove that $k\leq\fine(H^0_\mm(R/I))+1$, or, equivalently, that if $a$ is an integer such that $I_b=(I\sat)_b$ for all $b\geq a$
then $I:l^a=I:l^\infty$. It suffices to show that $I:l^\infty\subseteq I:l^a$. If  $z\in I:l^\infty$ 
then $zl^a\in I:l^\infty$, which is $I\sat$ by Remark \ref{somehelp}. 
More precisely $zl^a\in (I\sat)_{\geq a}=I_{\geq a}\subseteq I$ and therefore $z\in I:l^a$, as desired.
\end{proof}

\noindent  
In the next theorem we shall work with an ideal $I\subset K[X_1,\ldots,X_n]$ of height less than $n$, 
since this is the first non-trivial case. In fact, if $\dim R/I=0$ and $I$ is generated in degree $\leq d$, 
then $I$ contains a complete intersection of forms of degree at most $d$, therefore $\reg(I)\leq n(d-1)+1$.\\
We adopt the standard agreement that a product over the empty set is $1$.

\begin{theorem}\label{main}
Let $I$ be an homogeneous ideal of $K[X_1,\ldots,X_n]$ of height $c<n$ and  generated in degree  $\leq d$. 
Then, if  $l_n,\ldots,l_{c+1}$ is an almost-regular sequence of linear forms, one has
$$\reg(I)\leq\max\{d, \reg(I+(l_n))\}+d^c\prod_{i=c+2}^{n}\reg(I+(l_n,\ldots,l_{i})).$$
\end{theorem}

\begin{proof}
The proof of the theorem consists essentially in proving two separate inequalities:\\
\begin{equation}\label{A}
\reg(I)\leq \max\{d, \reg(I+(l_n))\} + \lambda\left(\frac{I:l_n}{I}\right);
\end{equation}
\begin{equation}\label{B}
\begin{split} 
\hbox{for}& \hbox{ all } i\geq c+2 \\
&\lambda\left(\frac{(I+(l_n,\ldots,l_{i+1})):l_i}{I+(l_n,\ldots,l_{i+1})}\right)\leq
\lambda\left(\frac{(I+(l_n,\ldots,l_i))\sat+(l_{i-1})}{I+(l_n,\ldots,l_{i-1})}\right)K_{i-1},
\end{split}
\end{equation}
where for all $i=1,\ldots,n$ the integer $K_i$ denotes the saturation index of $I+(l_n,\ldots,l_{i+1})$ with respect to $l_i$ i.e. the 
smallest integer $k$ such that $(I+(l_n,\ldots,l_{i+1})):l_i^\infty=(I+(l_n,\ldots,l_{i+1})):l_i^k$.\\

\noindent
{\em Proof of \eqref{A}:}\\
For the sake of notational simplicity we set $r\doteq\max\{d, \reg(I+(l_n)\}$ and $\lambda\doteq\lambda(I:l_n/I)$. We want to prove that
$I$ is $(r+\lambda)$-regular and for this purpose we make use of a well-known result of Bayer and Stillman (see for instance \cite{BaM} 
Theorem 3.3). By virtue of this criterion, if one wants to prove that an ideal $J$ generated in degrees $\leq d$ 
is $m$-regular, one has to check that $\left(\frac{J:l}{J}\right)_m=0$
and $m\geq\max\{d,\reg(J+(l))\}$, where $l$ is an
almost-regular linear form for $R/J$.
Thus,  we only have to prove that $(I:l_n/I)_{r+\lambda}=0$. Consider the chain 
$$\left(I:l_n/I\right)_{\geq r}\supset\left(I:l_n/I\right)_{\geq r+1}\supset\ldots\supset
\left(I:l_n/I\right)_{\geq r+\lambda}\supset\left(I:l_n/I\right)_{\geq r+\lambda+1}$$ 
and observe that if one of the inclusion is not strict, then one has that $(I:l_n/I)_{r+i}=0$ for some integer
$0\leq i\leq \lambda$, and since $r+i\geq r$ the criterion implies that $I$ is $(r+i)$-regular, therefore $(r+\lambda)$-regular. 
If this is would not be the case, i.e. all of the above inclusions are strict, one would obtain that $\lambda\geq\lambda+1$ which is a 
contradiction.\\

\noindent
{\em Proof of \eqref{B}:}\\
We first prove that, for all $i\geq c+1$, one has
\begin{equation}\label{B1}
\lambda\left(\frac{(I+(l_n,\ldots,l_{i+1})):l_i}{I+(l_n,\ldots,l_{i+1})}\right)=\lambda\left(\frac{(I+(l_n,\ldots,l_{i+1}))\sat+(l_i)}
{I+(l_n,\ldots,l_i)}\right).
\end{equation}
For simplicity fix an integer $c+1\leq i\leq n$ and let $J\doteq I+(l_n,\ldots,l_{i+1})$. Many applications of  Lemma \ref{snakesnake} 
yield 
\begin{equation*}
\begin{split}
\lambda\left(\frac{J:l_i}{J}\right) &=\lambda\left(\frac{J:l_i+(l_i)}{J+(l_i)}\right)+\lambda\left(\frac{J:l_i^2}{J:l_i}\right)\\
        &=\lambda\left(\frac{J:l_i+(l_i)}{J+(l_i)}\right)+\lambda\left(\frac{J:l_i^2+(l_i)}{J:l_i+(l_i)}\right)+
\lambda\left(\frac{J:l_i^3}{J:l_i^2}\right)\\
        &=\lambda\left(\frac{J:l_i+(l_i)}{J+(l_i)}\right)+\lambda\left(\frac{J:l_i^2+(l_i)}{J:l_i+(l_i)}\right)+\ldots+
        \lambda\left(\frac{J:l_i^{K_i}+(l_i)}{J:l_i^{K_i-1}+(l_i)}\right)\\
\end{split}
\end{equation*}
and since the sum is telescopic, by virtue of Remark \ref{somehelp}, one gets $\lambda\left(\frac{J:l_i}{J}\right)=
\lambda\left(\frac{J\sat+(l_i)}{J+(l_i)}\right)$, as required.\\
For all $i\geq c+2$ we now prove that 
\begin{equation}\label{B2}
\lambda\left(\frac{(I+(l_n,\ldots,l_{i+1}))\sat+(l_i)}{I+(l_n,\ldots,l_i)}\right)\leq\lambda\left(\frac{(I+(l_n,\ldots,l_i))\sat
+(l_{i-1})}{I+(l_n,\ldots,l_{i-1})}\right)K_{i-1}.
\end{equation}
For this purpose we want to study the length of a composition series of the module on the left-hand side of \eqref{B2}.
Since $(I+(l_n,\ldots,l_{i+1}))\sat+(l_i)$  is contained in 
$(I+(l_n,\ldots,l_i))\sat$ and $I+(l_n,\ldots,l_i)\subseteq (I+(l_n,\ldots,l_i)):l_{i-1}^{K_{i-1}}=(I+(l_n,\ldots,l_i))\sat$ 
(see again Remark \ref{somehelp}), 
we may consider the chain of inclusions 
$$I+(l_n,\ldots,l_i)\subset (I+(l_n,\ldots,l_i)):l_{i-1}\subset \ldots\subset (I+(l_n,\ldots,l_i)):l_{i-1}^{K_{i-1}}$$ in order to find 
such an estimate.\\ 
In the above chain we have exactly $K_{i-1}$ inclusions, thus
\eqref{B2} is proven if we are able to show that for all positive integers $a$ and $i\geq c+2$ 
one has:
\begin{equation}\label{B3}
\lambda\left(\frac{(I+(l_n,\ldots,l_i)):l_{i-1}^a}{(I+(l_n,\ldots,l_i)):l_{i-1}^{a-1}}\right)\leq 
\lambda\left(\frac{(I+(l_n,\ldots,l_i))\sat+(l_{i-1})}{I+(l_n,\ldots,l_{i-1})}\right).
\end{equation}
But the last is yielded by an iterated use of Lemma \ref{snakesnake} applied to the ideal $I+(l_n,\ldots,l_i)$. 
The last observation completes the proof of \eqref{B}.\\

\noindent
We now conclude the proof of the theorem. 
By \eqref{B1}, we can apply \eqref{B} iteratively  and obtain 
\begin{equation*}
\begin{split}
\lambda\left(\frac{I:l_n}{I}\right)\leq&\,\lambda\left(\frac{(I+(l_n,\ldots,l_{c+2}))\sat+(l_{c+1})}{I+(l_n,
\ldots,l_{c+1})}\right)K_{c+1}\cdot K_{c+2}
\cdot\ldots\cdot K_{n-1}\\
\leq & \,d^c\prod_{i=c+1}^{n-1}K_i.
\end{split}
\end{equation*}
The last inequality is due to the fact that $\lambda\left(\frac{(I+(l_n,\ldots,l_{c+2}))\sat+(l_{c+1})}{I+(l_n,
\ldots,l_{c+1})}\right)$ is less than or equal to that of $R/(I+(l_n,\ldots,l_{c+1}))$, which is an Artinian algebra. 
Therefore its length is bounded by 
$\lambda(K[X_1,\ldots,X_c]/(f_1,\ldots,f_c))$, where $f_1,\ldots,f_c$ is a regular sequence of $c$ elements of degree $\leq d$; thus 
the latter is limited above by $d^c$.\\
By virtue of \eqref{A}  we now  have $\reg(I)\leq \max\{d,\reg(I+(l_n))\} + d^c\prod_{i=c+1}^{n-1}K_i$ and it is sufficient to prove that
$K_i\leq\reg(I+(l_n,\ldots,l_{i+1}))$ for all $i=c+1,\ldots,n-1$, but this is an application of  Lemma \ref{scianti}.
\end{proof}

\begin{remark}\label{bia} Let  $I\subset K[X_1,\dots,X_n]$ be an homogeneous ideal generated in
degree $\leq d$. If the height of $I$ is $n$, we know already from what we discussed before Theorem \ref{main} that 
$\reg(I)\leq n(d-1)+1$. Furthermore, if $\Ht{I}=1$ then there exists an homogeneous polynomial 
$f$ of degree $0< a\leq d$ such that $I=(f)J$ and $J$ is an ideal generated in degree $\leq d-a$. 
Thus, the ideal $I$ is a shifted copy of $J$ and $\reg(I)=\reg(J)+a$.
\end{remark}
\noindent We are now in a  position of proving the sought after bounds.

\begin{corollary}\label{main2}
Let $I\subset K[X_1,\ldots,X_n]$  be an ideal of  height $c<n$ and generated in degree $\leq d$.  
Then  $$\reg(I)\leq (d^c+(d-1)c+1)^{2^{n-c-1}}.$$ 
\end{corollary}

\begin{proof} 
Let $l_n,\dots,l_{c+1}$ be an almost-regular sequence of linear forms.
By virtue of Theorem \ref{main} we are able to compute a bound for
the regularity of $(I+(l_n,\dots,l_i)$, $i\geq c+1$,  in the following way.
First we observe that the regularity of $I+(l_n,\ldots,l_i)$ equals
that of its image $\bar{I}$ in $K[X_1,\ldots,X_i]$ by restriction.
Moreover, the quotient algebra $R/(I+(l_n,\ldots,l_{c+1}))\simeq K[X_1,\ldots,X_c]/\bar{I}$ is Artinian  
and its  regularity is bounded by $c(d-1)+1$. We set $B_0\doteq (d-1)c+1$. Now we apply  Theorem
\ref{main} to the image of $I+(l_n,\dots,l_{c+2})$ in $K[X_1,\dots,X_{c+1}]$ and we obtain that
$\reg (I+(l_n,\dots,l_{c+2}))\leq (d-1)c+1+d^c$. We set  the latter to be $B_1$.
For all $i\geq 2$ we define recursively $B_i$ to be $B_{i-1}+\prod_{j=1}^{i-1}B_j$. It is easy to deduce
that $B_i=(B_{i-1}-B_{i-2})B_{i-1}+B_{i-1}\leq (B_{i-1})^2$.  Hence
$B_i\leq (B_1)^{2^{i-1}}$ for all $i\geq 1$ and $\reg(I)\leq 
B_{n-c} \leq ((d-1)c+1+d^c)^{2^{n-c-1}}$, as desired.
\end{proof}

\begin{corollary}
Let $I\subset K[X_1,\ldots,X_n]$  be an ideal generated in degree
$\leq d$. If $n=2$ then $\reg(I)\leq 2d-1$ otherwise, for $n\geq3$, we have
$$\reg(I)\leq ((d^2+2d-1)^{2^{n-3}}\leq (2d)^{2^{n-2}}.$$
\end{corollary}
\begin{proof}
The case $n=2$ is easy. If $n\geq 3$, we have only to verify that 
the worst possible situation occurs when the height of $I$ is $2$. Since the 
bounds are decreasing as a function of $c$, this is equivalent to saying that 
the case $\Ht(I)=1$ is not the worst possible, and this follows by what 
we discussed in Remark \ref{bia}. 
\end{proof}

\begin{example}
One could be interested in a slightly better estimate for the regularity and for this purpose could follow 
step-by-step the proof of Corollary \ref{main2}.\\
Consider for instance the case of the projective space $\PP^3$, i.e. the case $n=4$. As we said before, 
the worst possible case if provided by an ideal of height $2$. Since
$B_2=(B_1-B_0)B_{1}+B_{1}$, we have that the regularity of an homogeneous ideal in $\PP^3$ is bounded by
$((d^2+2d-1)-(2d-1))(d^2+2d-1)+(d^2+2d-1)=d^4+2d^3+2d-1$.
\end{example}

\noindent
{\bf Acknowledgments}\\
The authors would like to thank Aldo Conca for his helpful comments. 
The second author would also like to thank Markus Brodmann for having brought his attention to the subject.

\end{document}